\documentclass{siamart1116}

\usepackage{amsmath}
\usepackage{amsfonts}
\usepackage{amssymb}
\usepackage{graphicx, color}
\usepackage{hyperref}
\usepackage{epsfig}
\usepackage{epstopdf} 
\usepackage{comment}
\usepackage{cite}
\usepackage{enumerate}
\usepackage{siunitx} 

\definecolor{darkgreen}{rgb}{0,0.55,0}

\newtheorem{remark}{Remark}

\DeclareSymbolFont{AMSb}{U}{msb}{m}{n}
\DeclareMathSymbol{\N}{\mathbin}{AMSb}{"4E}
\DeclareMathSymbol{\Z}{\mathbin}{AMSb}{"5A}
\DeclareMathSymbol{\R}{\mathbin}{AMSb}{"52}
\DeclareMathSymbol{\Q}{\mathbin}{AMSb}{"51}
\DeclareMathSymbol{\I}{\mathbin}{AMSb}{"49}
\DeclareMathSymbol{\C}{\mathbin}{AMSb}{"43}

\newcommand{\verts}{\mathcal{H}(V)}
\newcommand{\edges}{\mathcal{H}(E)}
\newcommand{\vertso}{\mathcal{H}_{0}(V)}
\newcommand{\diverg}{\textnormal{div} }

\begin{document}

\title{Electrical Networks with Prescribed Current and Applications to Random Walks on Graphs}
\author{{Christina Knox\footnote{Department of Mathematics, University of California, Riverside, California, USA. E-mail: knox@math.ucr.edu.  }
\qquad Amir Moradifam \footnote{Department of Mathematics, University of California, Riverside, California, USA. E-mail: moradifam@math.ucr.edu. }}}
\date{}

\smallbreak \maketitle

\begin{abstract}
In this paper we study Current Density Impedance Imaging (CDII) on Electrical Networks. The inverse problem is to determine the conductivity matrix of an electrical network from the prescribed knowledge of the magnitude of the induced current along  the edges coupled with the imposed voltage or injected current on the boundary nodes. This problem leads to a weighted $l^1$ minimization problem for the corresponding voltage potential. We also investigate the problem of determining the transition probabilities of random walks on graphs from the prescribed expected net number of times the walker passes along the edges of the graph. Convergent numerical algorithms for solving such problems are also presented.  Our results can be utilized in the design of electrical networks when certain current flow on the network is desired as well as the design of random walk models on graphs when the expected net number of the times the walker passes along the edges is prescribed.  We also show that a mass preserving flow $J=(J_{ij})$ on a network can be uniquely recovered from the knowledge of $|J|=(|J_{ij}|)$ and the flux of the flow on the boundary nodes, where $J_{ij}$ is the flow from node $i$ to node $j$ and $J_{ij}=-J_{ji}$, and discuss its potential application in cryptography.

\end{abstract}
\maketitle

\section{Introduction} 
Let $G=(V,E)$ be a simple, undirected, weighted graph with $n$ vertices.  We can identify $G$ with an electrical network by placing a resistor with resistance $R_{ij}$ between every two vertices $i$ and $j$, for $0\leq i,j \leq n$ with $i\neq j$. We assign the weight $\sigma_{ij}=\frac{1}{R_{ij}}$ on each edge $E_{ij}$, and let $\sigma_{ij}=0$ if $i$ and $j$ are not connected. Suppose a voltage is applied to a subset of the vertices, denoted by $\partial V$ and called the boundary of $V$, then a current $J=(J_{ij})_{n \times n}$ will be induced on the edges of the graph,  where $J_{ij}$ is the current flowing from vertex $i$ to vertex $j$. In particular,  $J_{ij}=-J_{ji}$ and if the current flows from $i$ to $j$, then $J_{ij}>0$. We will also assume that $J_{ij}=0$ if the vertices $i$ and $j$ are not connected by an edge, and that $J_{ii}=0$. Note that $V=\partial V \cup \hbox{int}(V)=\{1,2,...,n\}$. We will view the voltage potential on $V$ as a vector $v=(v_1, v_2, ..., v_n) \in \mathbb{R}^n$ where $v_i$ is the voltage potential at vertex $i$.  We will also denote the imposed voltage potential on the boundary nodes by a function $f: \partial V \rightarrow \R$. By Kirchhoff's and Ohm's Law 
\begin{equation}
\label{eqn:kirkoff}
\sum_{j=1}^n \sigma_{ij}(v_i-v_j)=0 \ \ \hbox{for all } i \in \hbox{int}(V),
\end{equation}
where $\hbox{int}(V)=V \setminus \partial V$ are the interior nodes, and $v=f$ on $\partial V$ is the imposed voltage on the boundary nodes (Dirichlet boundary condition). Assume $((\sigma_{ij})_{n\times n}, f)$ is given on $E \times \partial V$.  Then \eqref{eqn:kirkoff} can be written as a system of $m=|\hbox{int}(V)|$ linear equations with $m$ unknowns, i.e. 
\begin{equation}
\label{FPmatrixA}
A_D v= b,
\end{equation}
where $v$ is a $m$ dimensional column vector containing the unknown voltage values at the interior nodes, $A_D$ is a $m \times m$ non-singular matrix (see Proposition \ref{ExUnofFP} below) depending on the conductivities, and $b$ is a $m$ dimensional column vector depending on the conductivities and the known voltage at the boundary. In particular the forward problem \eqref{eqn:kirkoff} always has a unique solution which is indeed the voltage potential associated to the conductivity problem on the network. 

On the other hand if a current $0 \neq g \in \R^{|\partial V|}$ is injected to the network on a subset of vertices $\partial V \subset V$ (Neumann boundary condition), then we necessarily have 
\begin{equation}
\sum_{i=1}^{|\partial V|} g_i=0, 
\end{equation}
and by Kirchhoff's and Ohm's Law the voltage potential $v$ satisfies 

\begin{eqnarray}\label{eqNeumann}
\left\{ \begin{array}{ll}
\sum_{j=1}^n \sigma_{ij}(v_i-v_j)=0 \ \ \ \  \hbox{for all } \ \  i \in \hbox{int}(V) \\\\
\sum_{j=1}^n \sigma_{ij}(v_i-v_j)=g_i\ \ \ \ \hbox{for all } \ \ i \in   \partial V.\\
\end{array} \right.
\end{eqnarray}
The above equations can be written as 
\begin{equation}
\label{FPmatrixAN}
A_N v= b,
\end{equation}
where $A_N$ is an $n\times n$ matrix depending on the conductivity $\sigma=(\sigma_{ij})_{n \times n}$, and $b$ is an $n$-dimensional column vector depending on the injected current on the boundary $\partial V$. The matrix $A_N$ also has unique solutions up to adding a constant (see Propositions \ref{Prop:A_N} and \ref{Prop:A_Nsoln} below) and the solution of \eqref{FPmatrixAN} is the voltage potential on the vertices of the graph. The matrix $A_N$ is in fact the well known graph laplacian of a weighted undirected graph.

As described above, the forward problems always have unique solutions up to a constant and can be easily solved by solving a linear system of equations. In this paper we are interested in the inverse problem of determining the conductivity matrix of an electrical network from the knowledge of the induced current along the edges of the network and Dirichlet or Neumann boundary conditions. This problem can also be understood as a design problem where one aims to design an electrical network that induces a prescribed current along its edges when a voltage $f\in \R^{|\partial V|}$ is applied to the boundary nodes $\partial V$, or when a current  $g\in \R^{|\partial V|}$ is injected on $\partial V$. These inverse problems are in the spirit of Current Density Imaging (CDI) and Current Density Impedance Imaging (CDII) in dimensions $n\geq 2$ which have been actively studied in recent years because of their potential applications in medical imaging, see \cite{Joy1, Joy2, HMN, JMN, Kim1, Kim2, Kim3, M-Indiana, MN, MNTim, MNTa_SIAM, NTT07, NTT08, NTT10, NTT11, NTV}. In dimension $n=3$ the induced current inside the conductive body $\Omega$ can be measured by Magnetic Resonance Imaging (MRI), see \cite{Joy1, Joy2}.

Random walks arise in many mathematical and physical models in biology, economics, computer and social networks, epidemiology, and statistical mechanics. Such models have been used to model infection on graphs such as spread of epidemics and rumours with mobile agents, see \cite{RW1, RW2}, voting patterns \cite{RWVoting,RWVoting2}, and stock market prices \cite{RWStocks}. Random walk models have also been proven to be a simple yet powerful method for extracting information from computer and social networks such as identification of reputable entities in a network. For instance Google's PageRank algorithm uses random walks to rank websites in their search engine results, see \cite{RW3, RW4}, and the survey papers \cite{Survey} and \cite{Sarkar}  for applications of random walks on graph in computer networks. Also see \cite{RWBook} for a wide variety of applications of random walks on graphs in statistical mechanics. The inverse problem we investigate here translate to intriguing questions in various contexts where a random walk model on graphs is utilized. The results could also be useful in the design of effective random walk models for achieving prescribed goals with random steps in a network. For instance, one can think of designing a random walk model with a prescribed high net number of times the walker passes along certain edges of the graph.

To the authors' best knowledge the natural inverse problem considered in this paper has not been studied elsewhere. In \cite{CurtMor} and \cite{ChungBeren}, the authors investigate the problem of recovering the conductivity of the edges from  the measurement of voltages at the boundary vertices, and measurements of the voltage, current, and conductivity on the boundary respectively. In \cite{CurtMor} the authors proved injectivity of this inverse problem for critical, circular and planar graphs and provided an explicit reconstruction method. Under the assumption of monotonicity of conductivities, partial uniqueness results are established in \cite{ChungBeren}. While the general theory of inverse problems on graphs is a rich field of study with applications in various disciplines, the above results are most closely related to this work.

There is a close connection between electrical networks and random walks on graphs (see \cite{DS}). In Section 5 we exploit this connection and apply our results on electrical networks to study the inverse problem of determining transition probabilities of random walk models from the net number of times the walker passes along the edges of the graph. We will also discuss a potential application of our results in public-key encryption, a seemingly unrelated  problem. 

The paper is organized as follows. In Section 2 we study the problem of determining the conductivity matrix of an electrical network from the knowledge of the magnitude of the induced current with Dirichlet boundary condition, and in Section 3 we study this problem with Neumann boundary data. In Section 4 we present a numerical algorithm for finding minimizers of the $l^1$ minimization problem we obtain in Sections 2 and 3.  In Section 5 the connection between random walks and electrical networks is discussed and we apply our results on electrical networks to the inverse problem of determining transition probabilities from the net number of time a random walker passes along the edges of the graph.

\section{Dirichlet Boundary Condition}
In this section we study the inverse problem of determining the conductivity matrix $\sigma=(\sigma_{ij})_{n \times n}$ from the knowledge of its induced current $J=(J_{ij})_{n \times n}$ on $E$ and the imposed voltage potential $f$ on $\partial V$ (Dirichlet  boundary conditions). Let $G=(V,E)$ be an undirected, simple, connected graph with $n$ vertices, and suppose a voltage is applied to some subset of the vertices inducing the current $J=(J_{ij})_{n \times n}$ on $E$. Throughout the paper $|J|$ denotes the matrix $|J|:=(|J_{ij}|)_{n \times n}$, we will refer to $|J|$ as a measurement matrix.

We first show that the forward problem has a unique solution, i.e. $A_D$ is non-singular. One can find a proof in \cite{CurtMor} and we present a brief proof for the sake of completeness. 
\begin{proposition}
\label{ExUnofFP}
The matrix $A_D$ is non-singular.
\end{proposition}
{\bf Proof.} For every $i \in int(V)$ it follows from \eqref{eqn:kirkoff} that $v_i$ is the weighted average of the voltage potential in its neighboring nodes, i.e.
\begin{equation}
v_i=\frac{\sum_{j=1}^n \sigma_{ij}v_j}{\sum_{j=1}^n \sigma_{ij}}.
\end{equation}
Consequently $v$ satisfies the strong maximum principle in the sense that if $v$ attains its maximum or minimum on an interior node, then $v$ must be constant on $V$. In particular, $v$ attains its minimum and maximum on the boundary $\partial V$. 

Now suppose  $A_D v=A_D\tilde{v}=b$. Then $w=v-\tilde{v}$ satisfies  
\[\sum_{j=1}^n \sigma_{ij}(w_i-w_j)=0 \ \ \hbox{for all } i \in \hbox{int}(V).\]
Since $w=0$ on $\partial V$, it follows from the above maximum principle that $w=0$ on $V$. Thus  the matrix $A_D$ is non-singular. \hfill $\square$ \\ \\
An immediate consequence of Proposition \ref{ExUnofFP} is that the forward problem \eqref{eqn:kirkoff} always has a unique solution.

\begin{definition}
We say that a vertex $i$ is an interior vertex and write $i \in \hbox{int}(V)$ if 
\[ J_i := \sum_{j=1}^{n} J_{ij}=0.\]
Otherwise we say that $i$ is boundary vertex and write $i \in \partial V$.  For every $i\in \partial V$, $J_i$ is the current flowing in ($J_i<0$) or out ($J_j>0$) of the graph at vertex $i$. In particular, $V= int(V) \cup \partial V$ and $int(V) \cap \partial V=\emptyset$. 
\end{definition}

\begin{definition}
\label{def:validpotential}
Given $f:\partial V \rightarrow \R$ and a measurement matrix $a=(a_{ij})_{n \times n}$ with $a_{ij} \in [0, \infty)$ for all $1\leq i,j \leq n$ and $a_{ij}=0$ when $i=j$ and $E_{i,j} \not\in E$, we say that a symmetric matrix $\sigma=(\sigma_{ij})_{n \times n}$ with $\sigma_{ij}\in [0,\infty]$ is a conductivity matrix associated to the data $(f,a)$, if there exists a function $v: \{1,2,..., n\}\rightarrow \R$ with $v|_{\partial V}=f$, and a matrix $J=(J_{ij})_{n \times n}$ such that 
\[ J_{ij}=\sigma_{ij}(v_i-v_j) \ \ \hbox{and} \ \ |J_{ij}|=a_{ij} \hbox{ for all } i,j \hbox{ with } v_i\neq v_j,\]
and 
\begin{equation*}
\sum_{j=1}^{n}J_{ij}=0
\end{equation*}
for all $i \in int(V)$. When $a_{ij}\neq 0$ and $v_i=v_j$, then we formally define $\sigma_{ij}=\infty$ and say that the edge between nodes $i$ and $j$ is a perfect conductor. We shall also refer to the function $v$ as a voltage potential and denote the set of all voltage potentials corresponding to the data $(f,a)$ by $\mathcal{V}_{(f,a)}$.  \\ 
\end{definition}

For any measurement matrix $a=(a_{ij})_{n \times n}$, define the function $I: \R^n\rightarrow \R$ by 
\begin{equation}
\label{eqn:functional}
I(u)=\frac{1}{2}\sum_{i,j}a_{ij} |u_i-u_j|,
\end{equation}
and for $f \in \R^{|\partial V|}$ consider the minimization problem 
\begin{equation}
\label{LGP}
\min \{ I(u):  \ \ u\in \R^n \ \ \hbox{and} \ \ u|_{\partial V}=f\}.
\end{equation}
We shall prove that $u \in \mathcal{V}_{(f,a)}$ if and only if it is a minimizer of the least gradient problem. Let us first study the dual of the minimization problem above. 

\subsection{The Dual problem}
Here we discuss the dual of the least gradient problem \eqref{LGP} and study the connection between these two problems.

Let $\mathcal{H}(V)$ be the set of all real valued functions on the vertices. We shall view a function $u \in \mathcal{H}(V)$ as a vector in $\mathbb{R}^n$. Also let $\mathcal{H}(E)$ to be the space of all functions on $E$, i.e. the space of all $n \times n$ matrices $b=(b_{ij})$, where $b_{ij}$ denotes the value of the function on the edge from vertex $i$ to $j$, with the additional convention that $b_{ij}=0$ if the edge from $i$ to $j$ is not in $E$, and  $b_{ii}=0$. 
\begin{definition}
Let $u,v \in  \verts$ and $a,b \in \edges$. Then we define the inner products 
\begin{align}
\langle u,v \rangle_{\verts}  = \sum_{i=1}^n u_iv_i, \hspace{3em} \langle b,d \rangle_{\edges}  = \sum_{i,j} b_{ij}d_{ij}
\end{align}
on $\verts \times \verts$ and $\edges \times \edges$, respectively. The spaces $\mathcal{H}(V)$ and $\mathcal{H}(E)$ equipped with the above inner products are Hilbert spaces.
\end{definition}
Next we define two linear operators $D : \verts \rightarrow \edges$ and $\diverg : \edges \rightarrow \verts$ which play crucial roles in our arguments. 
\begin{definition}
For $u \in \verts$ we define $D u \in \mathcal{H}(E)$ as
\begin{align}
(D u)_{ij}=u_i-u_j
\end{align}
if the edge connecting $i$ to $j$ is in $E$, and $0$ otherwise. Also for $b \in \edges$ we define $\diverg \hspace{1pt} b \in \verts$ as follows 
\begin{align}
(\diverg \hspace{1pt} b)_i=\sum_{j} b_{ji}-b_{ij}.
\end{align}
\end{definition}
Observe that if $b \in \edges$ is anti-symmetric, that is $b_{ij}=-b_{ij}$ for all $1 \leq i,j \leq n$, then the divergence is simply $-2\sum_{j}b_{ij}$. We shall refer to $D$ and $\diverg$ operators as gradient and divergence, respectively, since they play the role in our setting of the standard gradient and divergence operators on $\R^n$, $n\geq 2$. Note that the definition of the gradient and divergence given here does not depend on the weights (conductivities) of the graph as it would normally when defining these operators on a weighted graph. Since in the inverse problems we consider in this paper, the conductivities are unknown, these definitions are desirable.  Let us first show that $-\diverg$ is the adjoint of $D$. 
\begin{proposition}
Let $u \in \verts$ and $b \in \edges$. Then 
\begin{equation*}
\langle u, -\diverg \hspace{1pt} b \rangle _{\verts}=\langle  Du,b \rangle _{\edges}.
\end{equation*}
\end{proposition}
{\bf Proof.}  Let $u \in \verts$ and $b \in \edges$. Then
\begin{eqnarray*}
\langle u, -\diverg  \hspace{1pt} b \rangle_{\verts} &=& \sum_{i}u_i(-(\diverg b)_i)\\
&=& \sum_{i}u_i\sum_{j}(b_{ij}-b_{ji})\\
&=& \sum_{i} \sum_{j} u_ib_{ij}-\sum_{j}\sum_{i}u_j b_{ij}\\
&=& \sum_{i,j} (u_i-u_j)b_{ij}\\
&=& \sum_{i,j} (D u)_{ij}b_{ij}\\
&=&\langle D u,b \rangle _{\edges}.
\end{eqnarray*}
\hfill $\Box$

Let $f \in \R^{|\partial V|}$ and define 
\[\mathcal{H}_f=\{u\in \mathcal{H}(V): u|_{\partial V}=f  \}.\]
For $a \in \mathcal{H}(E)$ we take $a \geq 0$ to mean that every entry is non-negative. Then for $0 \leq a \in \mathcal{H}(E)$ and $f \in \R^{|\partial V|}$, the least gradient problem \eqref{LGP}  can be written as 
\begin{equation}\label{LGPP}
\underset{u \in \mathcal{H}_f}{\textnormal{min}} \hspace{2pt}\frac{1}{2}\sum_{i,j }a_{ij}|u_i-u_j|=\underset{u \in \mathcal{H}_f}{\textnormal{min}} \hspace{2pt} \frac{1}{2} \langle a, |D u| \rangle_{\edges},
\end{equation}
where we have used the notation $|Du|_{ij}=|(Du)_{ij}|$. Now choose $u_f \in \mathcal{H}_f.$ Define $\vertso \subset \mathcal{H}(V)$ to be the space of functions on $V$ which are equal to zero on $\partial V$. Then we can equivalently write the primal problem \eqref{LGPP} as
\begin{equation}\label{PPrime}
\underset{u \in \vertso}{\textnormal{min}} \hspace{2pt}\frac{1}{2}\sum_{i,j}a_{ij}|u_i-u_j+(u_f)_i-(u_f)_j|= \underset{u \in \vertso}{\textnormal{min}} \hspace{2pt} \frac{1}{2} \langle a, |D u+D u_f| \rangle_{\edges}.
\end{equation}
Define $F: \edges \rightarrow \mathbb{R}$ and $G: \vertso \rightarrow \mathbb{R}$ as follows 
\begin{equation}
F(d)=\frac{1}{2} \langle a,|d+D u_f| \rangle_{\edges} \ \ \hbox{and} \ \ G(u)\equiv 0.
\end{equation}
Then \eqref{PPrime} can be written as 
\begin{equation*}
(P) \ \ \ \ \alpha_P:=\underset{u \in \vertso}{\textnormal{min}} \hspace{2pt}F(D u)+G(u).
\end{equation*}
By Rockafellar-Fenchel duality (see \cite{ET}), this problem admits a dual problem which can be expressed as
\begin{equation}\label{dual}
\underset{b \in \mathcal{H}(E)}{\max}\hspace{2pt} -G^*(-\diverg b)-F^*(-b),
\end{equation}
where $F^*$ and $G^*$ denote the convex conjugate of $F$ and $G$, respectively.  It is easy to see that 
\begin{eqnarray*}
G^*(u) &=& \sup_{v\in \vertso} \sum_{i} u_iv_i\\
&=& \begin{cases} 0 & \textnormal{if} \hspace{5pt} u \equiv 0 \hspace{5pt} \textnormal{on} \hspace{5pt} int(V) \\
\infty &\textnormal{otherwise}.
\end{cases}
\end{eqnarray*}
Next we compute the convex conjugate of $F$. 
\begin{lemma}\label{F*}
Let $a=(a_{ij}) \in \edges$ with $a_{ij} \geq 0$ and $u_f \in \mathcal{H}_f(V)$. Then
\begin{equation}
F^*(b)=\begin{cases}-\langle b,D u_f \rangle_{\edges} & \textnormal{if} \hspace{1em} |b| \leq \frac{1}{2}a\\ \infty & \hbox{otherwise}.
\end{cases}
\end{equation}
\end{lemma}
{\bf Proof.} Suppose $|b| \leq \frac{1}{2}a$, that is $|b_{ij}| \leq \frac{1}{2}a_{ij}$ for all $i,j$. Then
\begin{eqnarray*}
F^*(b) &=& \sup_{d \in \mathcal{H}(E)} ( \langle d,b \rangle_{\edges}-\frac{1}{2}\langle a, |d+D u_f | \rangle_{\edges})\\
&=& - \langle b,D u_f \rangle_{\edges}+\sup_{d \in \mathcal{H}_a(E)} ( \langle d,b \rangle_{\edges}-\frac{1}{2}\langle a,|d| \rangle_{\edges})\\
&=& - \langle b,D u_f \rangle_{\edges}+\sup_{d \in \mathcal{H}_a(E)}( \sum_{i,j}d_{ij}b_{ij}-\frac{1}{2}a_{ij}|d_{ij}|)\\
&\leq & - \langle b,D u_f \rangle_{\edges}+\sup_{d \in \mathcal{H}_a(E)} \sum_{i,j}|d_{ij}|(|b_{ij}|-\frac{1}{2}a_{ij})\\
& \leq & - \langle b,D u_f \rangle_{\edges}.
\end{eqnarray*}
Taking $d=0$ we also get $F^*(b) \geq - \langle b,D u_f \rangle_{\edges}$.\\

Now suppose that there exists  $1\leq i_0,j_0 \leq n$ such that $|b_{i_0 j_0}| > \frac{1}{2}a_{i_0 j_0}$. Let $d_{i_0 j_0}=\lambda b_{i_0 j_0}$, and $d_{i j}=0$ otherwise, where $\lambda \in \R$. Then we have 
\begin{eqnarray*}
F^*(b) &=& - \langle b,D u_f \rangle_{\edges}+\sup_{d \in \mathcal{H}_a(E)} (\sum_{i,j}d_{ij}b_{ij}-\frac{1}{2}a_{ij}|d_{ij}|)\\
& \geq &- \langle b,D u_f \rangle_{\edges}+\sup_{\lambda >0} \lambda(b^2_{i_0 j_0}-\frac{1}{2}a_{i_0j_0}|b_{i_0j_0}|)\\
& = & - \langle b,D u_f \rangle_{\edges}+\sup_{\lambda >0} \lambda |b_{i_0 j_0}| (|b_{i_0j_0}|-\frac{1}{2}a_{i_0j_0})\\
&=& \infty. 
\end{eqnarray*}
\hfill $\Box$

Thus  the dual problem \eqref{dual} can be written as 
\begin{eqnarray*}
(D) \ \ \ \ \alpha_D:=\textnormal{sup}\{-\langle b,D u_f \rangle_{\edges} : \ \  b \in \edges,  \ \ |b| \leq \frac{1}{2}a,  \textnormal{ and } \hbox{div}(b)\equiv 0 \textnormal{ on } int(V)\}.
\end{eqnarray*}

Given that $u_i=0$ for at least one $i \in V$ one can show that any minimizing sequence of the the primal problem is uniformly bounded. Hence a convergent subsequence exists and a minimizer of the primal problem (P) always exists. On the other hand, it follows from Theorem III.4.1 in \cite{ET} that the dual problem (D) also has a solution. Indeed since  $I(u)=\frac{1}{2} \langle a,|Du+D u_f| \rangle_{\edges} $ is convex and $J: \mathcal{H}(E) \rightarrow \R$ with $J(p)=\frac{1}{2} \langle a,|p| \rangle_{\edges}$ is continuous at $p=0$, the condition (4.8) in the statement of Theorem III.4.1 in \cite{ET} is satisfied. The weighted $l^1$ minimization problem \eqref{LGP} does not have an unique minimizer and thus the conductivity inducing the current $J$ on $E$ is not unique. However we can characterize the non-uniqueness.

\begin{theorem}\label{duality}
The infimum of the primal problem (P) is equal to the supremum of the dual problem (D). Moreover, the dual problem has an optimal solution $b$, and $J=-2b$ satisfies 
\begin{equation}\label{currentvalue} |J_{ij}|=a_{ij} \hbox{ for every } i,j \hbox{ with } v_i\neq v_j
\end{equation}
and
\begin{equation}\label{currentdirect}
J_{ij}(v_i-v_j) \geq 0 \hbox{ for all }  1\leq i,j \leq n,
\end{equation}
for every minimizer $v$  of \eqref{LGP}. 
Conversely, if $u \in \mathcal{H}_f$ and the above equation holds then then $u$ is a minimizer of \eqref{LGP}.
\end{theorem}
{\bf Proof.} A solution $b$ to the dual problem always exists and the infimum of the primal problem (P) is equal to the supremum of the dual problem by Theorem III.4.1 in \cite{ET} as discussed above.
Let $v$ be a minimizer of \eqref{LGP}. Then 
\begin{eqnarray}\label{computations1}
\alpha_P= I(v)=\frac{1}{2}\sum_{i,j}a_{ij} |v_i-v_j|&\geq& \sum_{i,j}|b_{ij}| |v_i-v_j| \geq  \sum_{i,j}-b_{ij} (v_i-v_j)  \\
&=& \langle -b,Dv \rangle_{\edges} = \langle \diverg b, v \rangle_{\verts} \nonumber \\
&=&  \sum_{i \in \partial V}(\diverg b)_i v_i= \sum_{i \in \partial V}(\diverg b)_i f_i = \alpha_D=\alpha_P. \nonumber
\end{eqnarray}
Hence the inequalities in \ref{computations1} are indeed equalities and thus
\[|b_{ij}|=\frac{1}{2}a_{ij} \textnormal{ for every } i,j \textnormal{ with } v_i\neq v_j\]
and 
\[b_{ij}(v_i-v_j) \leq 0 \textnormal{ for all }  1\leq i,j \leq n.\]
Therefore if we let $J=-2b$ we we see that \eqref{currentvalue} and \eqref{currentdirect} hold. It is not hard to see that the converse also holds from the above computations. \hfill $\Box$

\begin{corollary}\label{cor1}
If $u$ and $v$ are two arbitrary minimizers of (\ref{LGP}), then
\[ (u_i-u_j)(v_i-v_j) \geq 0 \textnormal{ for all } 1\leq i,j \leq n. \]
\end{corollary}

\subsection{Voltage Potentials Have Minimum Energy}
We are now ready to prove the following theorem. 

 \begin{theorem}
\label{thm:iffvalidpotential}
Let $f$ be a function on $\partial V$ and $a$ be a measurement matrix. Then $v \in \mathcal{V}_{(f,a)}$ if and only if it is a minimizer of the least gradient problem \eqref{LGP}.
\end{theorem}
\textbf{Proof.} Suppose $v \in \mathcal{V}_{(f,a)}$ and let $J$ be the corresponding current on $E$. Then 
\begin{eqnarray}\label{computations}
\label{eqn:forwardmin}
I(v)=\frac{1}{2}\sum_{i,j}a_{ij} |v_i-v_j|&=&\frac{1}{2}\sum_{i,j}|J_{ij}| |v_i-v_j| \geq \frac{1}{2} \sum_{i,j}J_{ij} (v_i-v_j)  \\
&=& \sum _{i=1}^{n} v_i \sum_{j=1}^{n} J_{ij}=\sum _{i \in int (V)} v_i J_i+\sum _{i \in \partial V} v_i J_i \nonumber\\
&=&  \sum _{i \in \partial V} v_i J_i =  \sum _{i \in \partial V} f_i J_i. \nonumber
\end{eqnarray}
Therefore the minimum of the least gradient problem \eqref{LGP} is equal to $\sum _{i \in \partial V} f_i J_i$. Moreover the minimum is achieved for every  $v\in \mathcal{V}_{(f,|J|)}$. 

Now suppose $v$ is a minimizer of the problem \eqref{LGP} and let $b$ be a solution of the dual problem (D) and let $J=-2b$. Then by Theorem \ref{duality}
\[|J_{ij}|=a_{ij} \textnormal{ for all } i,j \textnormal{ with } v_i \neq v_j\]
and since $\diverg J=0$ on $int(V)$
\[\sum_{j=1}^n J_{ij}=0 \textnormal{ for all } i \in int(V). \]
For $v_i \neq v_j$ define $\sigma_{ij}=\frac{J_{ij}}{v_i-v_j}\geq 0$. Then 
\[J_{ij}=\sigma_{ij} (v_i-v_j) \ \ \ \ \hbox{for all} \ \ i,j \ \ \hbox{with}\ \ v_i \neq v_j.\] Thus $v \in \mathcal{V}_{(f,a)}$ and the proof is complete. \hfill $\Box$


\begin{remark}\label{remarkD}
Note that every minimizer $v$ of \eqref{LGP} uniquely determines a conductivity matrix $\sigma$. Corollary \ref{cor1} indicates that the directions of the flow of the current along the edges is unique, despite multiplicity of the minimizer of \eqref{LGP}. Indeed if two conductivity matrices $\sigma^1$ and $\sigma^2$ with $0\leq \sigma^1_{ij}, \sigma^2_{ij}< \infty$ induce the currents $J^1$ and $J^2$ on a network when the voltage $f$ is imposed on $\partial V$, and  $|J^1|=|J^2|$, then $J^1=J^2$. This is a counter-intuitive result.  
\end{remark}

\subsection{Multiple Measurements}
Suppose we have two data sets $(f^1,a^1)$ and $(f^2,a^2)$, and would like to find a conductivity matrix $\sigma$ inducing the currents with magnitudes $a^1$ and $a^2$, when the voltage potentials $f^1$ and $f^2$ are imposed on the boundary vertices $\partial V^1$ and $\partial V^2$, respectively.

 Let $I^1$ and $I^2$ be defined by Equation \eqref{eqn:functional} for $a^1$ and $a^2$ respectively and for $u=(u^1,u^2) \in \mathbb{R}^n \times \mathbb{R}^n$ define
\begin{equation}
\Phi(u^1,u^2) = \sum_{\mathcal{C}^2} {\left\lvert \frac{u^1_i-u^1_j}{|J_{ij}^1|}-\frac{u^2_i-u^2_j}{|J_{ij}^2|} \right\rvert}^2,
\end{equation}
where 
\[\mathcal{C}^2=\{(i,j): \ \ 1\leq i,j \leq n \ \ \hbox{and}\ \ J^1_{ij},J^2_{ij}\neq 0\}.\]
Define
\begin{equation}
\label{eqn:MultiMeas}
\mathcal{F}(u^1,u^2)=I^1(u^1)+I^2(u^2)+\Phi(u^1,u^2)
\end{equation}
and
\[\mathcal{A}:=\{(u^1,u^2)\in \R^n \times \R^n: \ \ u^1=f^1 \ \ \hbox{on } \partial V^1 \ \ \hbox{and}\ \ \ u^2=f^2 \ \ \hbox{on } \partial V^2\}.\]
Now consider 
\begin{equation}\label{LGP2}
\inf_{(u^1,u^2)\in \mathcal{A} } \mathcal{F}(u^1, u^2).
\end{equation}
It is easy to see that \eqref{LGP2} always has a minimizer. 

\begin{theorem}\label{thm:MultiMeas}
Let $(u^1,u^2)$ be a minimizer of \eqref{LGP2}. 
\begin{enumerate}
\item If there exists a conductivity matrix $\sigma$ which induces the current $J^i$ with $|J^i|=a^i$ when the voltage potential $f^i$ is imposed on the boundary, denoted $\partial^i V$, $i=1,2$, then $\Phi(u^1,u^2)=0$. Moreover, 
\[\sigma_{ij}=\frac{a^1_{ij}}{|u^1_i-u^1_j|} \ \ \hbox{for all} \ \ i,j \ \ \hbox{with} \ \ u^1_i\neq u^1_j,\]
and
\[\sigma_{ij}=\frac{a^2_{ij}}{|u^2_i-u^2_j|} \ \ \hbox{for all} \ \ i,j \ \ \hbox{with} \ \ u^2_i\neq u^2_j.\]
\item  If there doesn't exist a conductivity matrix $\sigma$ inducing the current $J^i$ with $|J^i|=a^i$ when the voltage potential $f^i$ is imposed on the boundary noted $\partial V^i$, $i=1,2$, then $\Phi(u^1,u^2)>0$.
\end{enumerate}
\end{theorem}
\textbf{Proof.} (1) Suppose there exists a conductivity matrix $\sigma$ producing the data $(f^1,a^1)$ and $(f^2, a^2)$.  It follows directly from Theorem \ref{thm:iffvalidpotential} that the set of minimizers of \eqref{LGP2} is equal to $\mathcal{V}_{(f^1,a^1)} \times \mathcal{V}_{(f^2,a^2)}$. So the first statement follows.\\
(2) Suppose $\Phi(u^1,u^2)=0$. Then $u^1$ and $u^2$ minimize $I^1$ and $I^2$ over the appropriate spaces and so by Theorem \ref{thm:iffvalidpotential}, $u^1 \in \mathcal{V}_{(f^1,a^1)}$ and $u^2 \in \mathcal{V}_{(f^2,a^2)}$ and thus they each have corresponding conductivity matrices $\sigma^1$ and $\sigma^2$ that generate currents $J^1$ and $J^2$ respectively. However $\Phi(u^1,u^2)=0$ implies that these conductivities are in fact equal. \hfill $\Box$\\ \\
Now suppose a finite data set of measurements is given:
\[(f^1,a^1), (f^2,a^2), ..., (f^k,a^k), \ \ k\geq 2.\]
Define 
\[I^l=\frac{1}{2}\sum_{ij}a^l_{ij} |u_i-u_j|, \ \ 1\leq l \leq k,\]
and 
\[\Phi^k(u^1, u^2,..., u^k)=\sum_{l=2}^{k}  \sum_{\mathcal{B}^l} {\left\lvert \frac{u^1_i-u^1_j}{|J_{ij}^1|}-\frac{u^l_i-u^l_j}{|J_{ij}^l|} \right\rvert}^2,\]
where 
\[\mathcal{C}^{l}=\{(i,j): \ \ 1\leq i,j \leq n \ \ \hbox{and}\ \ J^1_{ij},J^l_{ij}\neq 0\}.\]
Consider the weighted $l^1$ minimization problem 
\begin{equation}\label{LGPk}
\inf_{(u^1,u^2,..., u^k)\in \mathcal{A}^k} \sum_{l=1}^{k} I^{l}(v^l)+\Phi^k(u^1, u^2,..., u^k),
\end{equation}
where 
\[\mathcal{A}^k:= \{(u^1,u^2,..., u^k): \ \ u^{l}\in \R^n \ \ \hbox{and}\ \ u^l=f^l \ \ \hbox{on} \ \ \partial V^l, \ \ i=1,2,..., k\}.\]
One can similarly prove the following theorem. 

\begin{theorem}
Let $(u^1,u^2,..., u^k)$ be a minimizer of \eqref{LGPk}. 
\begin{enumerate}
\item If there exists a conductivity matrix $\sigma$ which induces the current $J^l$ with $|J^l|=a^l$ when the voltage potential $f^l$ is imposed on the boundary noted $\partial V^l$, $l=1,2, ...,k$, then $\Phi(u^1,u^2, ...,u^k)=0$. Moreover, 
\[\sigma_{ij}=\frac{a^l_{ij}}{|u^l_i-u^l_j|} \ \ \hbox{for all} \ \ i,j \ \ \hbox{with} \ \ u^l_i\neq u^l_j, \ \ l=1,2,...,k.\]
\item  If there doesn't exist a conductivity matrix $\sigma$ inducing the current $J^l$ with $|J^l|=a^l$ when the voltage potential $f^l$ is imposed on the boundary noted $\partial V^l$, $l=1,2,3,...,k$, then $\Phi(u^1,u^2, ..., u^k)>0$. 
\end{enumerate}
\end{theorem}

\section{Neumann Boundary Condition}
Let $G=(V,E)$ be an undirected simple connected graph with $n$ vertices, and suppose the current $0 \neq g\in \R^{|\partial V|}$ is injected to a subset $\partial V$ of $V$, regarded as boundary of $V$, inducing the current $J=(J_{ij})$ on $E$. Then $g$ should satisfy the compatibility assumption  
\begin{equation}\label{compatablity}
\sum_{i=1}^{|\partial V|}g_i=0.
\end{equation}
We will again denote $|J|:=(|J_{ij}|)_{n \times n}$ and refer to $|J|$ as a measurement matrix.  The following proposition characterizes solutions of the forward problem \eqref{eqNeumann}. \\ 
\begin{proposition}
\label{Prop:A_N}
Let $A_N$ be the matrix defined in \eqref{FPmatrixAN}. Then 
\[Ker (A_N)=\{(c,c, ..., c) \in \R^n: \ \ c\in \R\}.\]
\end{proposition}
{\bf Proof.} Suppose $A_N w=0$ for some $w \in \R^{n}$.  Then it follows from \eqref{eqNeumann} that
\begin{eqnarray*}
\frac{1}{2}\sum_{i,j}\sigma_{ij} (w_i-w_j)^2& = & \frac{1}{2} \sum _{i=1}^{n} w_i \sum_{j=1}^{n} \sigma_{ij}(w_i-w_j) -  \frac{1}{2} \sum _{j=1}^{n} w_j \sum_{i=1}^{n}  \sigma_{ij}(w_i-w_j) \\ 
&=&\sum _{i=1}^{n} w_i \sum_{j=1}^{n} \sigma_{ij}(w_i-w_j)  \\ 
&=& 0. 
\end{eqnarray*}
Hence $w_i=w_j$ for all $i$ and $j$ connected by an edge. Since $G$ is connected the proof is complete. \hfill $\Box$ \\ 
\begin{proposition}
\label{Prop:A_Nsoln}
The equation $A_Nv=b$ has a solution if and only if $\sum_{i=1}^n b_i=0$.
\end{proposition}
{\bf Proof.} By the Fredholm Alternative from linear algebra, $A_Nv=b$ has a solution if and only if $b \in Ker({A_N}^T)^{\perp}$. By the previous proposition and the fact that $A_N$ is symmetric we have
\[Ker({A_N}^T)^{\perp}=Ker(A_N)^{\perp}=\{b \in \mathbb{R}^N: \ \ \sum_{i=1}^n b_i=0\}.\]
\hfill $\Box$
 
Therefore if $\sum_{i=1}^n b_i=0$, up to adding a constant the equation \eqref{eqNeumann} has a unique solution. The following is the analog to Definition \ref{def:validpotential}.

\begin{definition}
\label{def:validpotentialN}
Given $0 \neq g:\partial V \rightarrow \R$ satisfying $\sum_{i=1}^{|\partial V|}g_i=0$ and a measurement matrix $a=(a_{ij})_{n \times n}$ with $a_{ij} \in [0, \infty)$ for all $1\leq i,j \leq n$ and $a_{ij}=0$ when $i=j$ and $E_{ij} \not\in E$, we say that a symmetric matrix $\sigma=(\sigma_{ij})_{n \times n}$ with $\sigma_{ij}\in [0,\infty]$ is a conductivity matrix associated to the data $(g,a)$, if there exists a function $v: \{1,2,..., n\}\rightarrow \R$ with and a matrix $J=(J_{ij})_{n \times n}$ such that 
\[ J_{ij}=\sigma_{ij}(v_i-v_j) \ \ \hbox{and} \ \ |J_{ij}|=a_{ij} \hbox{ for all } i,j \hbox{ with } v_i\neq v_j,\]
\[\sum_{j=1}^{n}J_{ij}=g_i \hbox{ for all } i \in \partial V\]
and
\begin{equation*}
\sum_{j=1}^{n}J_{ij}=0 \hbox{ for all } i \in int(V).
\end{equation*}
When $a_{ij}\neq 0$ and $v_i=v_j$, then we formally define $\sigma_{ij}=\infty$ and say that the edge between nodes $i$ and $j$ is a perfect conductor. We shall also refer to the function $v$ as a voltage potential and denote the set of all voltage potentials corresponding to the data $(g,a)$ by $\mathcal{V}_{(g,a)}$.  \\ 
\end{definition}

For a measurement matrix $a=(a_{ij})_{n \times n}$, define the function $I: \R^n\rightarrow \R$ by 
\begin{equation}
\label{eqn:functionalNeumann}
I(u)=\frac{1}{2}\sum_{i,j}a_{ij} |u_i-u_j|.
\end{equation}
Also for $g \in \R^{|\partial V|}$ satisfying \eqref{compatablity} define 
\[\mathcal{M}_g:=\{u \in \R^n: \sum_{i\in \partial V} u_i g_i=1\}.\]
We shall prove that the voltage potential is a minimizer of the $l^1$ minimization problem 
\begin{equation}\label{LGPNeumann}
\min_{u\in \mathcal{M}_g}\frac{1}{2}\sum_{i,j}a_{ij} |u_i-u_j|. 
\end{equation}
Let us first study the dual of this problem.  

\subsection{The Dual problem}
In this section we discuss the dual of the least gradient problem \eqref{LGPNeumann} and study its connection to the primal problem. Let $0 \neq g \in \R^{|\partial V|}$  satisfying \eqref{compatablity}. Choose $u_g \in \verts$ such that 
\[\sum_{i\in \partial V} (u_g)_i g_i=1.\]
Define 
\[\mathcal{M}_0:= \{u\in \mathcal{H}(V): \ \ \sum_{i\in \partial V} u_i g_i=0\}.\]
Then we can equivalently write the primal problem \eqref{LGPNeumann} as
\begin{equation}\label{LGPPrimeNeumann}
\underset{u \in \mathcal{M}_0}{\textnormal{min}} \hspace{2pt}\frac{1}{2}\sum_{i,j}a_{ij}|u_i-u_j+(u_g)_i-(u_g)_j|= \underset{u \in \mathcal{M}_0}{\textnormal{min}} \hspace{2pt} \frac{1}{2} \langle a, |Du+Du_g| \rangle_{\edges}.
\end{equation}
Define $F: \mathcal{H}(E) \rightarrow \mathbb{R}$ and $G: \mathcal{M}_0 \rightarrow \mathbb{R}$ as follows 
\begin{equation}
F(d)=\frac{1}{2} \langle a,|Du+Du_g| \rangle_{\edges} \ \ \hbox{and} \ \ G(u)\equiv 0.
\end{equation}
Then \eqref{LGPPrimeNeumann} can be written as 
\begin{equation*}
(P_N) \ \ \ \ \alpha_{P_N}:=\underset{u \in \mathcal{M}_0}{\textnormal{min}} \hspace{2pt}F(Du)+G(u).
\end{equation*}
As before this problem admits a dual problem which can be expressed as
\begin{equation}\label{dualNeumann}
\underset{b \in \mathcal{H}(E)}{\max}\hspace{2pt} -G^*(-\diverg b)-F^*(-b).
\end{equation}
From Lemma \ref{F*} we have 
\begin{equation*}
F^*(b)=\begin{cases}-\langle b,Du_g \rangle_{\edges} & \textnormal{if} \hspace{1em} |b| \leq \frac{1}{2}a\\ \infty & \textnormal{otherwise}.
\end{cases}
\end{equation*}
Next we compute $G^*$. 
\begin{lemma}
Let $G: \mathcal{M}_0 \rightarrow \mathbb{R}$ be defined as $G\equiv 0$. Then for $G^*: ( \mathcal{M}_0)^* \rightarrow \R$ we have

\begin{equation}G^*(D^* b)=\left\{ \begin{array}{ll}
0 &\hbox{if} \ \ b \in \mathcal{B}\\
\infty   &\hbox{otherwise},
\end{array} \right.
\end{equation}\label{G*} 
where 
\[\mathcal{B}:=\{b \in \mathcal{H}(E):  \ \  \diverg b\equiv 0 \hbox{ on } int(V) \hbox{ and } (\diverg b)_i=\lambda g_i \hbox{ for all } i \in \partial V, \hbox{ for some } \lambda \in \R\}.\]

\end{lemma}
\textbf{Proof.} First note that 
\begin{eqnarray*}
G^*(D^*b) &=& \sup_{u\in \mathcal{M}_0} \langle D^*b,u \rangle_{\verts} = \sup_{u\in \mathcal{M}_0} \langle b, Du \rangle_{\edges} = \sup_{u\in \mathcal{M}_0} -\langle \diverg b, u \rangle_{\verts} \\
&=& \begin{cases} 0 & \textnormal{if} \hspace{5pt} \diverg b \in {\mathcal{M}_0}^{\perp}\\
\infty &\textnormal{otherwise}.
\end{cases}
\end{eqnarray*}
Let  $h \in \verts$ with $h_i=0$ if $i \in int(V)$ and $h_i=g_i$ if $i \in \partial V$, and 
\[N=\{\lambda h: \lambda \in \R \} \subset \verts.\] 
Observe that $\mathcal{M}_0=\{u\in \verts : \langle h,u \rangle_{\verts}=0 \}$. Hence $\mathcal{M}_0=N^{\perp}.$ Since ${N^{\perp}}^{\perp}=N$, see \cite{Hal},
\[M_0^{\perp}=N, \]
and the result follows. \hfill $\square$ \\ \\
Therefore the dual problem \eqref{dualNeumann} can be written as 
\begin{eqnarray*}
(D_N) \ \ \ \ \alpha_{D_N}:=\sup_{b\in \mathcal{D}}\{-\langle b,Du_g \rangle_{\edges}\},
\end{eqnarray*}
where $\mathcal{D}=\{b \in \mathcal{B}: |b|\leq \frac{1}{2}a\}.$

Similar to before one can show that \eqref{LGPNeumann} has a minimizer. Similar to the Dirichlet boundary condition case, it follows from Theorem III.4.1 in \cite{ET} that the dual problem ($D_N$) also has a solution and characterizes the non-uniqueness of solutions of the primal problem \eqref{LGPNeumann}. 

\begin{theorem}\label{dualityN}
The infimum of the primal problem $(P_N)$ is equal to the supremum of the dual problem $(D_N)$. Moreover, the dual problem has an optimal solution $b$, and $J=-2b$ satisfies 
\begin{equation}\label{currentvalueN} |J_{ij}|=a_{ij} \hbox{ for every } i,j \hbox{ with } u_i\neq u_j
\end{equation}
and
\begin{equation}\label{currentdirectN22}
J_{ij}(u_i-u_j) \geq 0 \hbox{ for all }  1\leq i,j \leq n,
\end{equation}
for every minimizer $u$ of \eqref{LGPNeumann}. Conversely, if \eqref{currentvalueN} and \eqref{currentdirectN22}  hold for some $\mathcal{M}_g$, then then $u$ is a minimizer of \eqref{LGPNeumann}.
\end{theorem}
{\bf Proof.} Let $b$ be a solution to the dual problem with corresponding $\lambda \in \R$. Suppose $u$ is a minimizer of \ref{LGPNeumann}. Then 
\begin{eqnarray}\label{computations2}
\alpha_{P_N}= I(u)=\frac{1}{2}\sum_{i,j}a_{ij} |u_i-u_j|&\geq& \sum_{i,j}|b_{ij}| |u_i-u_j| \geq  \sum_{i,j}-b_{ij} (u_i-u_j)  \\
&=& \langle -b,Du \rangle_{\edges} = \langle \diverg b, u \rangle_{\verts} \nonumber \\
&=&  \lambda \sum_{i \in \partial V}g_i u_i=\lambda = \alpha_{D_N}=\alpha_{P_N}. \nonumber
\end{eqnarray}
Thus the inequalities in \eqref{computations2} are indeed equalities and taking $J=-2b$ we we see that \eqref{currentvalueN} and \eqref{currentdirectN22} hold. It is easy to see from the above compuations that the converse also holds. \hfill $\Box$

\begin{corollary}\label{cor2}
If $u$ and $v$ are two arbitrary minimizers of (\ref{LGPNeumann}), then
\[ (u_i-u_j)(v_i-v_j) \geq 0 \hbox{ for all } 1\leq i,j \leq n. \]\\
\end{corollary}

\subsection{Voltage Potentials Have Minimum Energy}
We can now prove the analog to Theorem \ref{thm:iffvalidpotential}.

 \begin{theorem}
\label{thm:iffvalidpotentialN}
Let $g\neq 0$ be a function on $\partial V$ satisfying \ref{compatablity} and $a$ be a measurement matrix. If $v \in \mathcal{V}_{(g,a)}$, then $v$ is a minimizer of the least gradient problem \eqref{LGPNeumann}. Conversely, given any $a=(a_{i,j})$ with $a_{i,j}\geq 0$ and $g \in \R^{|\partial V|}$ satisfying \eqref{compatablity}, if $v$ is a minimizer of the least gradient problem \eqref{LGPNeumann}, then $v\in \mathcal{V}_{(\lambda g, a)}$ for some $\lambda>0$.

\end{theorem}
\textbf{Proof.} Suppose $v \in \mathcal{V}_{(g,a)}$ and let $J$ be the corresponding current on $E$. Following similar computations as in the proof of Theorem \ref{thm:iffvalidpotential} we have 
\begin{eqnarray}\label{computationsN}
\label{eqn:forwardminN}
I(v)=\frac{1}{2}\sum_{i,j}a_{ij} |v_i-v_j|&=&\frac{1}{2}\sum_{i,j}|J_{ij}| |v_i-v_j| \geq \frac{1}{2} \sum_{i,j}J_{ij} (v_i-v_j)  \\
&=& \sum_{i \in \partial V} v_i g_i =1.\nonumber
\end{eqnarray}
Therefore the minimum of the least gradient problem \eqref{LGPNeumann} is equal to 1. Moreover the minimum is achieved for every  $v \in \mathcal{V}_{(g,|J|)}$. 

Now suppose $v$ is a minimizer of the problem \eqref{LGPNeumann} and let $b$ be a solution of the dual problem $(D_N)$ with the corresponding $\lambda \in \R$.   Let $J=-2b$. Then by Theorem \ref{dualityN} we see that $v\in \mathcal{V}_{(\lambda g, a)}$.  \hfill $\Box$

\begin{remark}\label{remarkN}
Note that Corollary \ref{cor2} indicates that the direction of the flow of the current along the edges is unique, despite multiplicity of the minimizers of \eqref{LGP} (see also Remark \ref{remarkD}). 
\end{remark}

\subsection{Multiple Measurements}
Suppose we have two data sets $(g^1,a^1)$ and $(g^2,a^2)$, and would like to find a conductivity matrix $\sigma$ inducing the currents with magnitudes $|J^1|$ and $|J^2|$, when the currents $g^1$ and $g^2$ are injected on the boundary vertices $\partial^1V$ and $\partial^2 V$, respectively. 
We can consider the minimization problem
\begin{equation}\label{LGP2N}
\inf_{(v^1,v^2)\in \mathcal{K} } F(v^1, v^2).
\end{equation}
where $F$ is defined by \eqref{eqn:MultiMeas} and 
\[\mathcal{K}:=\{(v^1,v^2)\in \R^n \times \R^n: \ \ \sum_{j=1}^n v^1_i g^1_i=1 \ \ \hbox{on } \ \ \hbox{and}\ \ \ \sum_{j=1}^n v^2_i g^2_i=1 \}.\]
The analog to Theorem \ref{thm:MultiMeas} can be formulated and proved in this setting and we can also similarly extend to a finite number of measurements.

\section{Algorithms for finding minimizers}
In this section we present numerical algorithms for finding minimizers of the $l^1$ minimization problems discussed in Sections 3 and 4, yielding voltage potentials for Dirichlet or Neumann boundary conditions.  The primal problem $(P_D)$ and $(P_N)$ can be written as 

\begin{equation}\label{GeneralMinimization}
\min_{\{u\in H, d\in \mathcal{H}(E)\}}  F(d) \ \ \ \ \hbox{subject to} \ \ \ \ Du=d,
\end{equation}
where $H=\mathcal{H}_0(V)$ for the Dirichlet case and $H=\mathcal{M}_0$ for the Neumann boundary problem. This leads to the unconstrained problem
\begin{equation}
\min_{\{u\in H, d\in \mathcal{H}(E) \}} F(d)+\frac{\alpha}{2} \|Du-d\|^2. 
\end{equation}
To solve the above minimization problem, we use and develop an algorithm in the spirit of the alternating Split Bregman method which was first introduced by Goldstein and Osher \cite{GO}. The Split Bregman algorithm suggests initiating the vectors $b^0$ and $d^0$, and producing the sequences $u^k$, $b^k$, and $d^k$ as follows

\begin{eqnarray}\label{SB}
(u^{k+1},d^{k+1})&=&\hbox{argmin}_{u\in H, d\in H(E)} \{F(d)+\frac{\alpha}{2}\parallel b^k+Du -d\parallel^2_{2}\}, \\
b^{k+1}&=& b^k+Du^{k+1}-d^{k+1},\nonumber
\end{eqnarray}
where $\alpha>0$. Since the joint minimization problem (\ref{SB}) in both $u$ and $d$ is in general expensive to solve exactly, Goldstein and Osher \cite{GO} proposed the following 
Alternating Split Bregman algorithm for solving problems of type $(\ref{GeneralMinimization})$
\begin{eqnarray}
u^{k+1}&=&\hbox{argmin}_{u \in H} \parallel b^k+Du-d^k\parallel^{2}_2,\\
d^{k+1}&=&\hbox{argmin}_{d \in \mathcal{H}(E)} \{F(d)+\frac{\alpha}{2}\parallel b^k+Du^{k+1}-d\parallel^{2}_2\},\\
b^{k+1}&=&b^k+Du^{k+1}-d^{k+1}.
\end{eqnarray}
See \cite{Es, COS, GO, GM, Setzer, Setzer1} for more details. It is pointed out by Esser \cite{Es} and Setzer \cite{Setzer1} that the above idea to minimize alternatingly was first presented for the augmented Lagrangian algorithm by Gabay and Mercier \cite{GM} and Glowinski and Marroco \cite{GMar}. The resulting algorithm is called the alternating direction method of multipliers (ADMM) \cite{Gab} and is equivalent to the alternating split Bregman algorithm. The convergence of ADMM in finite dimensional Hilbert spaces was established  by Eckstein and Bertsekas \cite{EcBe}.  This in particular implies convergence of the alternating split Bregman algorithm in finite dimensional Hilbert spaces. Cai, Osher, and Shen \cite{COS} and  Setzer \cite{Setzer, Setzer1} also independently presented convergence results for the alternating split Bregman in finite dimensional Hilbert spaces.  In \cite{MN} and \cite{MNTim} the authors proved the convergence of the alternating split Bregman algorithm in infinite dimensional Hilbert spaces by showing that the alternating split bregman algorithim corresponds to the Douglas-Rachford splitting algorithm for the dual problem. Indeed the dual problems \eqref{dual} and \eqref{dualNeumann} can be written in the form 
\begin{equation}\label{inclusion}
0\in A(-b)+B(-b),
\end{equation}
where $A:=\partial G^* o(-\hbox{div})$ and $B=\partial F^*$ are maximal monotone operators on $H$. For a set valued operator $P: H\rightarrow 2^H$, let $J_P$ denote its resolvent, i.e. $J_P=(Id+P)^{-1}$.  Douglas-Rachford splitting algorithm states that for any initial elements $x_0$ and $p_0$ and any $\alpha>0$, the sequences $p_{k}$ and $x_k$ generated by the following algorithm
\begin{eqnarray}\label{DRS-18}
x_{k+1}&=&J_{\alpha A}(2p_k-x_k)+x_k-p_k \nonumber \\
p_{k+1}&=&J_{\alpha B}(x_{k+1}),
\end{eqnarray}
converges to some $x$ and $p$ respectively. Furthermore $p=J_{\alpha B}(x)$ and $p$ satisfies
\[0\in A(p)+B(p).\]
Let us introduce the sequences $b^k$ and $d^k$ with 
\[x_k=\alpha(b^k+d^k) \ \ \hbox{and}\ \ p_k=\alpha b_k.\]
Notice that both sequences $b^k$ and $d^k$ converge. The resolvents $J_{\alpha A}(2p_k-x_k)$ and $J_{\alpha B}(x_{k+1})$ can be computed as follows 
\begin{equation}
J_{\alpha A}(2p_k-x_k)=\alpha(b^k+Du^{k+1}-d^k)
\end{equation}
and
\begin{equation}
J_{\alpha B}(x_{k+1})=\alpha(b^k+Du^{k+1}-d^{k+1}),
\end{equation}
where $u^{k+1}$ and $d^{k+1}$ are minimizers of 
\begin{equation}
I_1(u)= \sum_{i,j} | b^k_{ij}+(D u)_{ij}-d^k_{ij} |^2
\end{equation}
and
\begin{equation}
I_2(d)=\frac{1}{2} \sum_{i,j} a_{ij}|d_{ij}+(D u_f)_{ij}| +\frac{\alpha}{2} \sum_{i,j} | b^k_{ij}+(D u^{k+1})_{ij}-d_{ij} |^2
\end{equation} 
over $u \in \vertso$ for the Dirichlet problem and over $u\in \mathcal{M}_0$ for the Neumann problem, and over $d \in \edges$.

In the case of Dirichlet boundary condition the minimizer of $I_1$ should satisfy the Euler-Lagrange equation 
\begin{eqnarray}\label{ELAlgorithmD}
\left\{ \begin{array}{ll}
\sum \limits_{j=1}^{n} (Du)_{ij} = \frac{1}{2}[(\diverg b^k)_i-(\diverg d^k)_i], \ \ \ \ \forall i \in \hbox{int}(V) \\
u_i = 0\ \ \ \ \hbox{for all } \ \ i \in \partial V.\\
\end{array} \right.
\end{eqnarray}
It follows from Proposition \ref{ExUnofFP} that the above system is uniquely solvable. 

In the case of Neumann boundary condition, $I_1$ also has a unique minimizer in $\mathcal{M}_0$ up to adding a constant, but identifying the solutions is more subtle.  First note that if $u$ is a minimizer $I_1$ in $\mathcal{M}_0$, then it satisfies the Euler-Lagrange equation 
\begin{eqnarray}\label{ELAlgorithmN}
\left\{ \begin{array}{ll}
\sum \limits_{j=1}^{n} (Du)_{ij} = \frac{1}{2}[(\diverg b^k)_i-(\diverg d^k)_i], \ \ \ \ \forall i \in \hbox{int}(V) \\
\sum\limits_{j=1}^{n}(Du)_{ij}=\beta g_i+[\frac{1}{2}(\diverg b^k)_i-(\diverg d^k)_i], \ \ \ \ \hbox{for all } \ \ i \in \partial V\\
\end{array} \right.
\end{eqnarray}
for some $\beta \in \R$. Conversely for $\beta \in \R$, every solution of the above equation which belongs to $\mathcal{M}_0$ is a minimizer of $I_1$.   Since $\sum_{i \in \partial V}g_i=0$ and $\sum_{i=1}^n (\diverg c)_i=0$ for any $c \in \edges$, by Propositions \ref{Prop:A_N} and \ref{Prop:A_Nsoln} the system \eqref{ELAlgorithmN} has a unique solution in $\mathcal{H}(V)$ for every $\beta \in \R$, up to adding a constant. To identify $\beta$ and find a solution of \eqref{ELAlgorithmN}  in $\mathcal{M}_0$, let $z$ be a solution of 

\begin{eqnarray}\label{fixedGap!}
\left\{ \begin{array}{ll}
\sum \limits_{j=1}^{n} (Dz)_{ij} =0, \ \ \ \ \forall i \in \hbox{int}(V) \\
\sum\limits_{j=1}^{n}(Dz)_{ij}=g_i \ \ \ \ \hbox{for all } \ \ i \in \partial V.\\
\end{array} \right.
\end{eqnarray}
Then 
\begin{eqnarray*}
0<\frac{1}{2}\sum_{i,j} (D_z)_{ij}& = & \frac{1}{2} \sum _{i=1}^{n} w_i \sum_{j=1}^{n} (Dz)_{ij} -  \frac{1}{2} \sum _{j=1}^{n} z_j \sum_{i=1}^{n}  (Dz)_{ij} \\ 
&=&\sum _{i=1}^{n} z_i \sum_{j=1}^{n}(Dz)_{i,j}  \\ 
&=&\sum _{i \in \partial V} z_i g_i.\\
\end{eqnarray*}
Hence
\[\sum _{i \in \partial V} z_i g_i>0.\]
Now let $u$ be a solution of 
\begin{eqnarray}
\sum \limits_{j=1}^{n} (Du)_{ij} = \frac{1}{2}[(\diverg b^k)_i-(\diverg d^k)_i], \ \ \ \ \forall i \in V. \
\end{eqnarray}
Define 
\[\beta:=-\frac{\sum _{i \in \partial V} u_i g_i}{\sum _{i \in \partial V} z_i g_i}.\]
Then $v=u+\beta z$ belongs to $\mathcal{M}_0$ and satisfies the equation \eqref{ELAlgorithmN}, and hence $v$ is the unique minimizer of $I_1$ over $\mathcal{M}_0$, up to adding a constant.

The minimizer of $I_2$ for the Dirichlet problem can be directly computed as
\begin{eqnarray}
d^{k+1}_{ij}=\left\{ \begin{array}{ll}
\max\{|w_{ij}|-\frac{a_{ij}}{2\alpha},0\}\frac{w_{ij}}{|w_{ij}|}-(D u_f)_{ij} \ \ \ \ \hbox{if } w_{ij}\neq 0 \\
-(D u_f)_{ij} \ \ \ \ \hbox{if } w_{ij}= 0, \\
\end{array} \right.
\end{eqnarray}
where $w_{ij}=(Du^{k+1})_{ij}+(D u_f)_{ij}+ b^k_{ij}$. For the Neumann problem $u_f$ is replaced by $v_g$.\\ \\
Therefore Douglas-Rachford splitting leads to the following convergent algorithms for the Dirichlet and Neumann problems. \\ \\
\textbf{Algorithm 1 (Finding a minimizer of the Dirichlet Problem)}\\ \\
Let $\alpha>0$, $u_f \in \verts$ with $u=f$ on $\partial V$ and initialize $b^0,d^0 \in \edges$. For $k\geq 0$:
\begin{enumerate}

\item Solve
\begin{eqnarray}
\left\{ \begin{array}{ll}
\sum\limits_{j} (Du^{k+1})_{ij} = \frac{1}{2}[(\diverg b^k)_i-(\diverg d^k)_i], \ \ \ \ \forall i \in \hbox{int}(V) \\
u^{k+1}_i = 0\ \ \ \ \hbox{for all } \ \ i \in \partial V.\\
\end{array} \right.
\end{eqnarray}

\item Compute $d^{k+1}$ 
\begin{eqnarray}
d^{k+1}_{ij}=\left\{ \begin{array}{ll}
\max\{|w_{ij}|-\frac{a_{ij}}{2\alpha},0\}\frac{w_{ij}}{|w_{ij}|}-(D u_f)_{ij} \ \ \ \ \hbox{if } w_{ij}\neq 0 \\
-(D u_f)_{ij} \ \ \ \ \hbox{if } w_{ij}= 0, \\
\end{array} \right.
\end{eqnarray}
where $w_{ij}=(D u^{k+1})_{ij}+(D u_f)_{ij}+ b^k_{ij}$. 
\item
Set
\[b^{k+1}_{ij}=b^k_{ij}+(D u^{k+1})_{ij}-d^{k+1}_{ij}.\]\\
\end{enumerate}
The following proposition follows directly from the convergence of Douglas-Rachford splitting algorithm and Theorem 1.2 in \cite{MN}. See also \cite{COS, Setzer, Setzer1}. 

\begin{proposition}\label{AlgorithmPropD}
Let $u^k$ $b^k$, and $d^k$ be the sequences produced by the Algorithm 1. Then $u^k\rightarrow u$ and $b^k\rightarrow\frac{1}{2 \alpha} J$, where $u$ and $J$ are solutions of the \eqref{PPrime} and it's dual problem (D), respectively. In addition $d^k\rightarrow Du$. In particular $u$ is a voltage potential corresponding to the data $(f,a)$ and $J$ is the induced current with $|J|=a$. \\ 
\end{proposition} 
\textbf{Algorithm 2 (Finding a minimizer of the Neumann Problem)}\\ \\
Let $\alpha>0$, $v_g \in \verts$ with $\sum \limits_{i \in \partial V}v_ig_i=1$ and initialize $b^0,d^0 \in \edges$. Also let $z\in \R^n$ be a solution of \eqref{fixedGap!} with $z_1=0$.   For $k\geq 0$:
\begin{enumerate}

\item (a) Solve
\begin{eqnarray}
\left\{ \begin{array}{ll}
\sum \limits_{j} (Du^{k+1})_{ij} = \frac{1}{2}[(\diverg b^k)_i-(\diverg d^k)_i], \ \ \ \ \forall i \in V \\
\end{array} \right.
\end{eqnarray}
with $u^{k+1}_1=0$.\\
(b) Compute 
\[\beta^{k+1}=-\frac{\sum _{i \in \partial V} u^{k+1}_i g_i}{\sum _{i \in \partial V} z_i g_i}\]
and set $v^{k+1}=u^{k+1}+\beta^{k+1} z$. 

\item Compute $d^{k+1}$ 
\begin{eqnarray}
d^{k+1}_{ij}=\left\{ \begin{array}{ll}
\max\{|w_{ij}|-\frac{a_{ij}}{2\alpha},0\}\frac{w_{ij}}{|w_{ij}|}-(D v_g)_{ij} \ \ \ \ \hbox{if } w_{ij}\neq 0 \\
-(D v_g)_{ij} \ \ \ \ \hbox{if } w_{ij}= 0, \\
\end{array} \right.
\end{eqnarray}
where $w_{ij}=(D v^{k+1})_{ij}+(D v_g)_{ij}+ b^k_{ij}$. 
\item
Set
\[b^{k+1}_{ij}=b^k_{ij}+(D v^{k+1})_{ij}-d^{k+1}_{ij}.\]\\ 
\end{enumerate}

Convergence of Douglas-Rachford splitting algorithm implies the following convergence result, see Theorem 1.2 in \cite{MN} and \cite{COS, Setzer, Setzer1}. 
\begin{proposition}\label{AlgorithmPropN}
Let $v^k$ $b^k$, and $d^k$ be the sequences produced by the Algorithm 2. Then $v^k\rightarrow v$ and $b^k\rightarrow \frac{1}{2\alpha}J$, where $v$ and $J$ are solutions of the \eqref{LGPPrimeNeumann} and it's dual problem ($D_N$), respectively. In addition $d^k\rightarrow Dv$. In particular $v$ is a voltage potential corresponding to the data $(\lambda g,a)$ for some $\lambda \in \mathbb{R}$ and $J$ is the induced current with $|J|=a$. Moreover $\lambda$ is the optimal values of the primal and dual problems $(P_N)$ and $(D_N)$, i.e. $\lambda=\alpha_{P_N}=\alpha_{D_N}$. \\ 
\end{proposition}
\subsection{Numerical Simulations}
We performed a set of numerical simulations in MATLAB to demonstrate convergence of Algorithm 1 and 2. A simple graph with 100 vertices was generated and edges were randomly assigned between nodes with a approximate density of 0.125. Random numbers uniformly distributed between 0 and 1 were then assigned to each edge as their conductivity. We then selected 5 boundary nodes and randomly assigned values between 0 and 1 as boundary data.  For the Dirichlet boundary data, the forward problem was solved to determine the current $J$, generating the data $a=|J|$. To generate the boundary data for the Neumann problem we found the current entering/leaving the system at each boundary vertex. The simulations for both the Dirichlet and Neumann boundary data were done on the same graph structure with the same current data  $|J|$. The nonsingular linear systems in algorithm 1 were solved using the MATLAB mldivide function and the singular linear systems in algorithm were solved using the pinv function. The vector $u_f$ was chosen to be zero on $int(V)$ and $f$ on the $\partial V$. The vector $v_g$ in Algorithm 2 was chosen using the MATLAB mldivide function. Tables 1 and 2 show the numerical errors for algorithms 1 and 2 on the same graph for different levels of tolerance. Simulations were run on a late 2013 MacBook Pro with a 2.4 GHZ Intel Core i5 processor. We used the $L^2$ matrix norm for error computations.

\begin{table}[h!]
  \begin{center}
    \caption{Numerical errors for algorithm 1 on 100 node graph with 1121 edges}
    \begin{tabular}{l|l|l|l}
      \textbf{Tolerance} & \textbf{Relative L2 Error} & \textbf{Number of Iterations} & \textbf{Elapsed Time(s)}\\ 
      \hline
      $10^{-3}$ & 1.2171$\times 10^{-3}$ &  16 & 0.069309\\ 
      $10^{-4}$ & 1.3160$\times 10^{-4}$& 22 & 0.102846\\ 
      $10^{-5}$ & 1.4494$\times 10^{-5}$ & 92 & 0.358250\\
      $10^{-6}$ & 1.3615$\times 10^{-6}$  & 133 & 0.405979 
    \end{tabular}
  \end{center}
\end{table}

\begin{table}[h!]
  \begin{center}
    \caption{Numerical errors for algorithm 2 on 100 node graph with 1121 edges}
    \begin{tabular}{l|l|l|l}
      \textbf{Tolerance} & \textbf{Relative L2 Error} & \textbf{Number of Iterations} & \textbf{Elapsed Time(s)}\\ 
      \hline
      $10^{-2}$ & 1.3069$\times 10^{-3}$ &   7  & 0.055400\\ 
      $10^{-3}$ & 1.3908$\times 10^{-4}$& 9& 0.071342\\ 
      $10^{-4}$ & 1.0235$\times 10^{-5}$ & 12 & 0.086956\\
      $10^{-5}$ & 1.1987$\times 10^{-6}$  & 24 & 0.147310 
    \end{tabular}
  \end{center}
\end{table}

While running our simulations we observed that the speed of convergence of Algorithm 1 varied quite wildly depending on the choice of boundary data. We also observed that the speed of convergence of Algorithm 2 was always the same or faster than that of Algorithm 1. To test this observation, we ran algorithms 1 and 2 on the same graph used in Tables 1 and 2 for 1000 different choices of Dirichlet boundary. The average number of iterations for each algorithm is shown in Table 3. We also remark that changing the structure of the graph also effects the speed of convergence. It is not clear to the authors that why Algorithm 2 converges faster than Algorithm 1, and an in depth analysis of the speed of convergences of algorithms 1 and 2 remain open.

\begin{table}[h!]
  \begin{center}
    \caption{Average Number of Iterations}
    \begin{tabular}{l|l|l|l}
      \textbf{Tolerance} & \textbf{Algorithm 1} & \textbf{Algorithm 2} \\ 
      \hline
      $10^{-3}$ & 21.175 &  15.918\\ 
      $10^{-4}$ & 46.097  & 18.905\\ 
      $10^{-5}$ & 111.847  & 23.486\\
      $10^{-6}$ & 227.624  & 32.846
    \end{tabular}
  \end{center}
\end{table}

\vspace{3cm}

\section{Applications}
In this section we discuss potential applications of our results on electrical networks on random walks on graphs and Cryptography. 
\subsection{Random Walks on Graphs}
Let $G=(V,E')$ be a connected, directed, and simple graph with $n$ nodes and consider a random walk on $G$. Suppose a random walker begins at node $a$ and walks until they reach node $b$ and if they return to $a$ before reaching $b$ they keep walking. Let $P=(P_{ij})\in \mathcal{H}(E)$ be the matrix of transition probabilities, i.e. $0 \leq P_{ij} \leq 1$ is the probability of the random walker walking from node $i$ to node $j$. In particular $\sum_{j}P_{ij}=1$ for all $1\leq i\leq n$. Let $W_{ij}$ be the expected number of times the walker walks from node $i$ to node $j$ before exiting the graph at node $b$. Note that $W_{ij}=-W_{ji}$. Can one determine transition probabilities $P=(P_{ij})$ from the knowledge of the boundary vertices $\{a,b\}$ and $W=(W_{i,j})$? In this section, among other results, we show that the answer is yes, and describe an algorithm for determining such $P$. 

There is a close connection between electrical networks and random walks on graphs \cite{DS}. Let $G=(V,E)$ be an electrical network with conductivity matrix $\sigma=(\sigma_{ij})$, $\sigma_{i,j}\in [0,\infty)$, and let $\partial V=\{a,b\}$.  Suppose a current $g$ with $g(a)=1$ and $g(b)=-1$ is injected to the network inducing a current $J$ along the edges. Define 
\begin{eqnarray}\label{ProbabilityConductivity}
\sigma_i:=\sum_{j=1}^{n}\sigma_{ij} \ \ \hbox{and}\ \  P_{ij}=\frac{\sigma_{ij}}{\sigma_i}
\end{eqnarray}
and assign the transition probability matrix $P$ to the graph $G=(V,E')$. Then the net number of times the walker taking an step from node $i$ to node $j$ is indeed $J_{ij}$, i.e. 
\[J=W.\]
Therefore if the boundary nodes $\partial V=\{a,b\}$ and the magnitude of expected net number of times the walker should walk along the edges of the graph is prescribed, by the method presented in Section 5, one can first find a conductivity matrix $\sigma$ inducing the current $J=W$ on network and compute transition probability matrix $P$ by \eqref{ProbabilityConductivity}.

The connection between random walks on graphs and electrical networks with Neumann boundary condition can be generalized to the case when $\partial V=\Gamma_a \cup \Gamma_b$ with $\Gamma_a \cap \Gamma_b= \emptyset$ and $\Gamma_a, \Gamma_b \neq \emptyset$. Let $g\in \R^{|\partial V|}$ with $g|_{\Gamma_a}\geq 0$ and $g|_{\Gamma_b}\leq 0$ and 
\[\sum_{i\in \Gamma_1}g_i=1 \ \ \hbox{and}\ \ \ \ \sum_{i\in \Gamma_2}g_i=-1.\]
Suppose we would like to determine a transition matrix $P$ such that if a random walker enters the network from a vertex $k$ in $\Gamma_a$ with probability $g_k$, then 
\begin{itemize}
\item  they exit the network at a node $l\in \Gamma_b$ with probability $|g_l|$
\item the expected net number of times they pass from vertex $i$ to node $j$ before exiting the network is $W_{ij}$, $1 \leq i,j \leq n$.
\end{itemize}
As explained above, to determine the transition matrix $P$ it suffices to find a conductivity matrix $\sigma$ inducing the current $J=W$  with Neumann data $g$ on $\partial V$. Then $P$ can be computed from \eqref{ProbabilityConductivity}. 

Suppose $\partial V=\{a,b\}$ and consider the inverse problems of determining the transition probabilities from the relative net number of times the walker walks between the edges of the graphs, i.e. $\alpha W=(\alpha W_{i,j})$ where $\alpha$ is a unknown constant. Then one can determine a transition probability $P$ by finding a conductivity matrix $\sigma$ by minimizing the $l^1$ minimization problem \eqref{eqn:functional} with $a=\alpha W$, $f(a)=1$ and $f(b)=0$. A transition matrix can also be obtained by minimizing \eqref{LGPNeumann} with the Neumann boundary condition $g(a)=1$ and $g(b)=-1$. 

\begin{remark}
Note that in this section we assume that the conductivity matrix $\sigma=(\sigma_{ij})$ satisfies $\sigma_{i,j}\in [0,\infty)$. Indeed we do not allow perfect conductors as otherwise the probability matrix $P$ in \eqref{ProbabilityConductivity} will not be well-defined. As described in the introduction, if for a minimizer $v$ of \eqref{LGP} or \eqref{LGPNeumann} we have $v_i=v_j$ and $|J_{i,j}|\neq 0$ for some $1\leq i,j\leq n$, then the edge $(i,j)$ is a perfect conductor, i.e. $\sigma_{i,j}=\infty$. If $v$ is minimizer of \eqref{LGP} or \eqref{LGPNeumann} leading to perfect conductance on an edge, then one may look for an increasing function $F: \R\rightarrow \R$ such that $u=(u_1,u_2,...,u_n):=(F(v_1),F(v_2),..., F(v_n))$ satisfies $u_i\neq u_j$ for $i\neq j$. Note that such $u$ will also be a minimizer of \eqref{LGP} or \eqref{LGPNeumann} and would provide a conductivity matrix $\sigma$ with $\sigma_{ij}\in [0,\infty)$, and hence the transition probabilities can be computed from \eqref{ProbabilityConductivity}. If such increasing function $F$ does not exists, then there exists no transition probability matrix $P$ for which the expected number of times the walker passes along the edges is $W$. 
\end{remark}

\subsection{Applications in Cryptography}
In this section we discuss a potential application of our results on electrical networks in public-key encryption. As stated in Remark \ref{remarkN}, Theorem \ref{thm:iffvalidpotential} implies that a mass preserving flow $J=(J_{ij})$ along the edges of a graph $G=(V,E)$ can be recovered from the knowledge of $|J|=(|J_{ij}|)$ and its net flux on the boundary nodes $\partial V$. More precisely, suppose $J_{i,j}$ is the current from node $i$ to node $j$ ($J_{ij}=-J_{ji}$ for $(i,j)\in E$), and suppose  
\[\sum_{j=1}^{n}J_{ij}=0 \ \ \ \ \hbox{for every interior node} \ \ i \not \in \partial V\]
and 
\[\sum_{j=1}^{n}J_{ij}=f_i \ \ \ \ \hbox{for every boundary node} \ \ i \in \partial V.\]
Then $J$ can be reconstructed from the knowledge of $(|J|,f, \partial V)$. This counter-intuitive result has a potential application in cryptography. To see the connection, let us translate a special case of this result to the language of matrices. 

Let $I_n$ be a subset of $ \{1,2, ..., 2n+1\}$ with $n$ elements and $\mathcal{A}_{I_n}$ be the space of $(2n+1) \times (2n+1)$ anti-symmetric matrices $A=(a_{ij})$ satisfying the following properties: \\

\begin{enumerate}[I]
\item. $a_{ij}\in \{-1,0,+1\}$ for $i\neq j$ and $a_{ii}=0$, for all $1\leq i,j \leq 2n+1$
\item. All rows of $A$ contain an even number of non-zero entries
\item. Sum of the entries of the $i$th row is equal to zero if $i \not\in I_n$
\item. For $i\in I_n$, the sum of the entries of the $i$th row is denoted by $f_i$, which is not necessarily zero. \\
\end{enumerate}

Note that $f \in \R^n$. Suppose a pair of communicators have agreed on a set of indices $I_n \subset \{1,2,...,2n+1\}$ with $n$ elements, both are aware of $I_n$, and would like to securely communicate a matrix $A\in \mathcal{A}_{I_n}$. Then the first party can just send the key $(|A|,f)$ where $f \in R^{n}$ is the sum of the entries of the rows of $A$ that belong to $I_n$. The second party can decrypt the message and find $A$ from the knowledge of $(|A|,f,I_n)$, using the algorithm we developed in Section 4. Since $a_{ij}$ only takes integer values in $\{-1,0,+1\}$, a few iterations of the algorithm should be enough to determine $A$. On the other hand, finding $A$ from the knowledge of $(|A|,f)$ would be extremely difficult for an adversary who is not aware of $I_n$. Indeed since all rows of $|A|$ have an even number entries equal to $1$, the adversary could not determine the boundary nodes $I_n$ from $|A|$. To decrypt the message, the adversary faces the problem of guessing  $I_n$ among ${ {2n+1}\choose{n}} $ subsets of $\{1,....,2n+1\}$ with $n$ elements and matching it with $f$. The number of different possibilities are 
\[n!{ {2n+1}\choose{n}} \simeq \frac{2^{2n+1}}{\sqrt{\pi n}}n!, \]
which grows very fast and makes the decryption for adversaries extremely difficult for large $n$. The above application in public-key encryption and the challenges of its implementation will be further studied in a forthcoming paper.

\bibliography{ElectNetwork}
\bibliographystyle{siamplain}
\end{document}